\documentclass{amsart}

\usepackage{amsfonts}
\usepackage{amsmath}
\usepackage{amsthm}
\usepackage{amssymb}
\usepackage[english]{babel}
\usepackage{graphics}
\usepackage{latexsym}
\usepackage{longtable} \usepackage{mathrsfs}
\usepackage{graphicx}
\usepackage{wrapfig} %\usepackage{txfonts}
\usepackage[pdftex,colorlinks=true,linkcolor=blue,citecolor=red]{hyperref}

\newtheorem{theorem}{Theorem}
\newtheorem{corollary}[theorem]{Corollary}
\newtheorem{lemma}[theorem]{Lemma}

\newtheorem{claim}[theorem]{Claim}
\newtheorem{example}[theorem]{Example}

\theoremstyle{definition}
\newtheorem{definition}[theorem]{Definition}

\newcommand{\mO}{\mathcal{O}}

\newcommand{\mD}{\mathcal{D}}

\newcommand{\A}{\mathrm{A}}
\newcommand{\R}{\mathbb{R}}
\newcommand{\N}{\mathbb{N}}

\newcommand{\X}{\textbf{X}}

\newcommand{\D}{\mathrm{D}}

\newcommand{\noi}{\noindent}
\newcommand{\ms}{\medskip}

\newcommand{\al}{\alpha}
\newcommand{\be}{\beta}

\newcommand{\de}{\delta}
\newcommand{\De}{\Delta}

\newcommand{\la}{\lambda}

\newcommand{\Om}{\Omega}

\newcommand{\larrow}{\longrightarrow}
\newcommand{\ot}{\otimes}

\newcommand{\p}{\partial}
\newcommand{\sub}{\subseteq}
\newcommand{\set}{\setminus}
\newcommand{\by}{\times}
\newcommand{\rk}{\mathrm{rk}}

\newcommand{\tr}{\mathrm{tr}}
\newcommand{\sgn}{\mathrm{sgn}}

\newcommand{\Div}{\mathrm{Div}}

\newcommand{\inter}{\mathrm{int}}

\newcommand{\bt}{\begin{theorem}}\newcommand{\et}{\end{theorem}}
\newcommand{\bd}{\begin{definition}}\newcommand{\ed}{\end{definition}}
\newcommand{\bl}{\begin{lemma}}\newcommand{\el}{\end{lemma}}
\newcommand{\beq}{\begin{equation}}\newcommand{\eeq}{\end{equation}}
\newcommand{\bc}{\begin{claim}}\newcommand{\ec}{\end{claim}}
\newcommand{\bex}{\begin{example}}\newcommand{\eex}{\end{example}}
\newcommand{\bcor}{\begin{corollary}}\newcommand{\ecor}{\end{corollary}}
\newcommand{\bp}{\begin{proof}}\newcommand{\ep}{\end{proof}}

\newcommand{\BPL}{\medskip \noindent \textbf{Proof of Lemma} }

\newcommand{\BPCOR}{\medskip \noindent \textbf{Proof of Corollary} }

\newcommand{\BPT}{\medskip \noindent \textbf{Proof of Theorem} }

\numberwithin{equation}{section}
\begin{document}

\title[Rigidity and flatness of maps with tangential Laplacian]{Rigidity and flatness of the image of certain classes of mappings having tangential Laplacian}

%\title[Rigidity and flatness of the image of $\infty$-Harmonic  maps]{Rigidity and %flatness of the image of certain classes of $\infty$-Harmonic and $p$-Harmonic maps}

\author{Hussien Abugirda, Birzhan Ayanbayev and Nikos Katzourakis}

\address[H. Abugirda]{Department of Mathematics , College of Science, University of Basra, Basra, Iraq AND Department of Mathematics and Statistics, University of Reading, Whiteknights, \\ PO Box 220, Reading RG6 6AX, UK}
%\curraddr{Department of Mathematics and Statistics, Case Western
%Reserve University, Cleveland, Ohio 43403}
\email{h.a.h.abugirda@pgr.reading.ac.uk}


\address[B. Ayanbayev]{Department of Mathematics and Statistics, University of Reading, Whiteknights, PO Box 220, Reading RG6 6AX, UK}
%\curraddr{Department of Mathematics and Statistics, Case Western
%Reserve University, Cleveland, Ohio 43403}
\email{ b.ayanbayev@pgr.reading.ac.uk}

%    Information for third author

\address[N. Katzourakis, corresponding author]{Department of Mathematics and Statistics, University of Reading, Whiteknights, PO Box 220, Reading RG6 6AX, UK}
%\curraddr{Department of Mathematics and Statistics, Case Western
%Reserve University, Cleveland, Ohio 43403}
\email{ n.katzourakis@reading.ac.uk}

  \thanks{\!\!\!\!\!\!\texttt{N.K. has been partially financially supported through the EPSRC grant EP/N017412/1}}

%    General info

%\date{}

\dedicatory{This paper is dedicated to Gunnar Aronsson with the utmost esteem for his pioneering work.}


\begin{abstract} In this paper we consider the PDE system of vanishing normal projection of the Laplacian for $C^2$ maps $u : \mathbb{R}^n \supseteq \Omega \longrightarrow \mathbb{R}^N$:
\[
[\![\mathrm{D} u]\!]^\bot \Delta u = 0 \ \, \text{ in }\Omega.  
\]
This system has discontinuous coefficients and geometrically expresses the fact that the Laplacian is a  vector field tangential to the image of the mapping. It arises as a constituent component of the $p$-Laplace system for all $p\in [2,\infty]$. For $p=\infty$, the $\infty$-Laplace system is the archetypal equation describing extrema of supremal functionals in vectorial Calculus of Variations in $L^\infty$. Herein we show that the image of a solution $u$ is piecewise affine if either the rank of $\mathrm{D} u$ is equal to one or $n=2$ and $u$ has the additively separated form $u(x,y)=f(x)+g(y)$. As a consequence we obtain corresponding flatness results for the images of $p$-Harmonic maps, $p\in [2,\infty]$.

\end{abstract}

\maketitle

\section{Introduction} \label{section1}

Suppose that $n,N$ are integers and $\Om$ an open subset of $ \R^n$. In this paper we study geometric aspects of the image $u(\Om)\sub \R^N$ of certain classes of $C^2$ vectorial solutions $u:\R^n \supseteq \Om\longrightarrow \R^N$ to the following nonlinear degenerate elliptic PDE system: 
\beq
\label{1.1}
[\![\D u]\!]^\bot \De u  = 0 \, \ \text{ in }\Om.
\eeq
Here, for the map $u$ with components $(u_1,...,u_N)^\top$ the notation $\D u$ symbolises the gradient matrix
\[
\ \ \ \ \D u(x) = \big(\D_i u_\al(x)\big)_{i=1...n}^{\al=1...N} \, \in\, \R^{N\by n}\ , \ \ \D_i \equiv \p /\p x_i,
\]
$\De u$ stands for the Laplacian
\[
\De u(x) = \sum_{i=1}^n \D^2_{ii}u(x) \, \in \, \R^N
\]
and for any $X\in \R^{N\by n}$, $ [\![X]\!]^\bot$ denotes the orthogonal projection on the orthogonal complement of the range of linear map $X :\R^n \larrow \R^N$:
\beq \label{1.2}
[\![X]\!]^\bot := \textrm{Proj}_{\mathrm{R}(X)^\bot}.
\eeq
Our general notation will be either self-explanatory, or otherwise standard as e.g.\ in \cite{E,D}. Note that, since the rank is a discontinuous function, the map $[\![\,\cdot\, ]\!]^\bot$ is discontinuous on $\R^{N\by n}$; therefore, the PDE system \eqref{1.1} has \emph{discontinuous coefficients}. The geometric meaning of \eqref{1.1} is that \emph{the Laplacian vector field $\De u $ is tangential to the image $u(\Om)$} and hence \eqref{1.1} is equivalent to the next statement: there exists a vector field
\[
\A\, : \, \R^n \supseteq \Om \longrightarrow \R^n
\]
such that 
\[
\Delta u  = \D u \,\A \, \ \ \text{ in }\Om. 
\]
As we show later, the vector field is generally discontinuous (Lemma \ref{lemma2}).

Our interest in \eqref{1.1} stems from the fact that it is a constituent component of the $p$-Laplace PDE system for all $p\in[2,\infty]$. Further, contrary perhaps to appearances, \eqref{1.1} is in itself a \emph{variational PDE system} but in a non-obvious way. Deferring temporarily the specifics of how exactly \eqref{1.1} arises and what is the variational principle associated with it, let us recall that, for $p\in[2,\infty)$, the celebrated $p$-Laplacian is the divergence system
\beq \label{1.33}
\ \ \ \De_p u:= \Div \big(|\D u|^{p-2}\D u \big) = 0\ \ \text{ in }\Om
\eeq
and comprises the Euler-Lagrange equation which describes extrema of the model $p$-Dirichlet integral functional
\beq  \label{1.34}
\ \ \ E_p(u) := \int_\Om|\D u|^p, \ \ \ u\in W^{1,p}(\Om,\R^N),
\eeq
in conventional vectorial Calculus of Variations. Above and subsequently, for any $X\in \R^{N\by n}$, the notation $|X|$ symbolises its Euclidean (Frobenius) norm: 
\[
|X| = \left(\sum_{\al=1}^N\sum_{i=1}^n\, (\X_{\al i})^2\right)^{\!\! 1/2}. 
\]
The pair \eqref{1.33}-\eqref{1.34} is of paramount important in applications and has been studied exhaustively. The extremal case of $p\to \infty$ in \eqref{1.33}-\eqref{1.34} is much more modern and intriguing, in that totally new phenomena arise which are not present in the scalar case. It turns out that one then obtains the following nondivergence PDE system
\beq 
\label{1.3}
\ \ \Delta_\infty u  := \Big(\D u \ot \D u + |\D u|^2[\![\D u]\!]^\bot \! \ot \mathrm{I}\Big):\D^{2}u = 0 \ \ \text{ in } \Om,
\eeq
which is known as the $\infty$-Laplacian. In index from, \eqref{1.3} reads
\[
\ \ \ \sum_{\be=1}^N\sum_{i,j=1}^n \Big(\D_i u_\al \, \D_j u_\be + |\D u|^2[\![\D u]\!]_{\al\be}^\bot\, \de_{ij}\Big)\D^2_{ij}u_\be\,=\,0, \  \ \ \al=1,...,N.
\]
The system \eqref{1.3} plays the role of the Euler-Lagrange equation and arises in connexion with variational problems for the supremal functional 
\beq \label{1.6}
\ \ \ E_\infty(u,\mathcal{O}) :=  \| \mathrm{D} u\|_{L^\infty(\mathcal{O})}, \ \ \ u \in W^{1,\infty}(\Om,\mathbb{R}^N), \ \mathcal{O} \Subset \Omega.
\eeq
The scalar case of $N=1$ in  \eqref{1.3}-\eqref{1.6} was pioneered by G.\ Aronsson in the 1960s \cite{A1}-\cite{A5} who initiated the field of Calculus of Variations in $L^\infty$, namely the study of supremal functionals and of their associated equations describing critical points. Since then, the field has developed tremendously and there is an extensive relevant literature (see e.g.\ \cite{BEJ, BJW1, BJW2, BL, BN, CDP, GM, P, P2} and the lecture notes \cite{B,C,K7}). In particular, although vectorial supremal  functionals began to be explored early enough, the $\infty$-Laplace system \eqref{1.3} which describes the necessary critical conditions in $L^\infty$ in the vectorial case $N \geq 2$ first arose in the early 2010s in \cite{K1}. The area is now developing very rapidly due to both the mathematical significance as well as the importance for applications in several areas (see \cite{ETT, AK, AyK, CKP, KP}, \cite{K2}-\cite{K4}, \cite{K8}-\cite{K10}).

In this paper we focus on the $C^2$ case and establish the geometric rigidity and flatness of the images of solutions $u:\R^n \supseteq \Om\longrightarrow \R^N$ to the nonlinear system \eqref{1.1}, under the assumption that either $\D u$ has rank at most 1, or that $n=2$ and $u$ has an additively separated form, see \eqref{asp}. As a consequence, we obtain corresponding flatness results for the images of solutions to \eqref{1.33} and  \eqref{1.3}. Both aforementioned classes of solutions furnish particular examples which provide substantial intuition for the behaviour of general extremal maps in Calculus of Variations in $L^\infty$, see e.g.\ \cite{A6,A7,C,K7,K3,K4,KP} where solutions of this form have been studied. Obtaining further information for the still largely mysterious behaviour of $\infty$-Harmonic maps is perhaps the greatest driving force to isolate and study the particular nonlinear system \eqref{1.1}. For example, it is not yet know to what extend the possible discontinuities of the coefficients relates to the failure of absolute minimality.

Let us note that the rank-one case includes the scalar and the one-dimensional case (i.e.\ when $\min\{n,N\}=1$), although in the case of $N=1$ (in which the single $\infty$-Laplacian reduces to $\D u \ot \D u :\D^2 u=0$) \eqref{1.1} has no bearing since it vanishes identically at any non-critical point. 

\emph{The effect of \eqref{1.1} to the flatness of the image can be seen through the $L^\infty$ variational principle introduced in \cite{K2}}, wherein it was shown that solutions to \eqref{1.1} of constant rank can be characterised as those having minimal area with respect to \eqref{1.6}-\eqref{1.34}. More precisely, therein the following result was proved: 

\begin{theorem}[cf.\ \text{\cite[Theorem 2.7, Lemma 2.2]{K2}}] Given $N\geq n\geq 1$, let $u : \R^n \supseteq \Om \larrow \R^N$ be a $C^2$ immersion defined on the open set $\Om$ (more generally $u$ can be a map with constant rank of its gradient on $\Om$). Then, the following statements are equivalent:

\begin{enumerate}
\item \label{1} The map $u$ solves the PDE system \eqref{1.1} on $\Om$.

\smallskip

\item \label{2} For all $p\in[2,\infty]$, for all compactly supported domains $\mO \Subset \Om$ and all $C^1$ vector fields $\nu : \mO \larrow \R^N$ which are normal to the image $u(\mO)\sub\R^N$ (without requiring to vanish on $\p\mO$), namely those for which $\nu = [\![\D u]\!]^\bot \nu$ in $\mO$, we have
\[
\|\D u \|_{L^p(\mO)}  \leq \, \|\D u + \D \nu \|_{L^p(\mO)}.
\]

\item \label{3} The same statement as in item \eqref{2} holds, but only for \emph{some} $p\in[2,\infty]$.
\end{enumerate}
If in addition $p<\infty$ in \eqref{2}-\eqref{3}, then we may further restrict the class of normal vector fields to those satisfying $\nu|_{\p \mO}=0$ (see Figure 1).
\end{theorem}

In the paper \cite{K2}, it was also shown that in the conformal class, \emph{\eqref{1.1} expresses the vanishing of the mean curvature vector of $u(\Om)$}.
 
The effect of \eqref{1.1} to the flatness of the image can be easily seen in the case of $n=1\leq N$ as follows: since
\[
\ \ [\![u']\!]^\bot u''  = 0 \ \, \text{ in }\Om \sub \R
\]
and in one dimension we have
\[
[\![u']\!]^\bot =\, \left\{
\begin{array}{ll}
\mathrm{I} -\dfrac{u'\ot u'}{|u'|^2}, & \text{ on }\{u'\neq0\},
\\
\mathrm{I}, & \text{ on }\{u' = 0\},
\end{array}
\right.
\]
we therefore infer that $u''=fu'$ on the open set $\{u'\neq0\} \sub \R$ for some function $f$, readily yielding after an integration that $u(\Om)$ is necessarily contained in a piecewise polygonal line of $\R^N$. As a generalisation of this fact, our first main result herein is the following:
\[
\underset{\text{\small \phantom{\Big|} Figure 1. Illustration of the variational principle characterising \eqref{1.1}.}}{\includegraphics[scale=0.17]{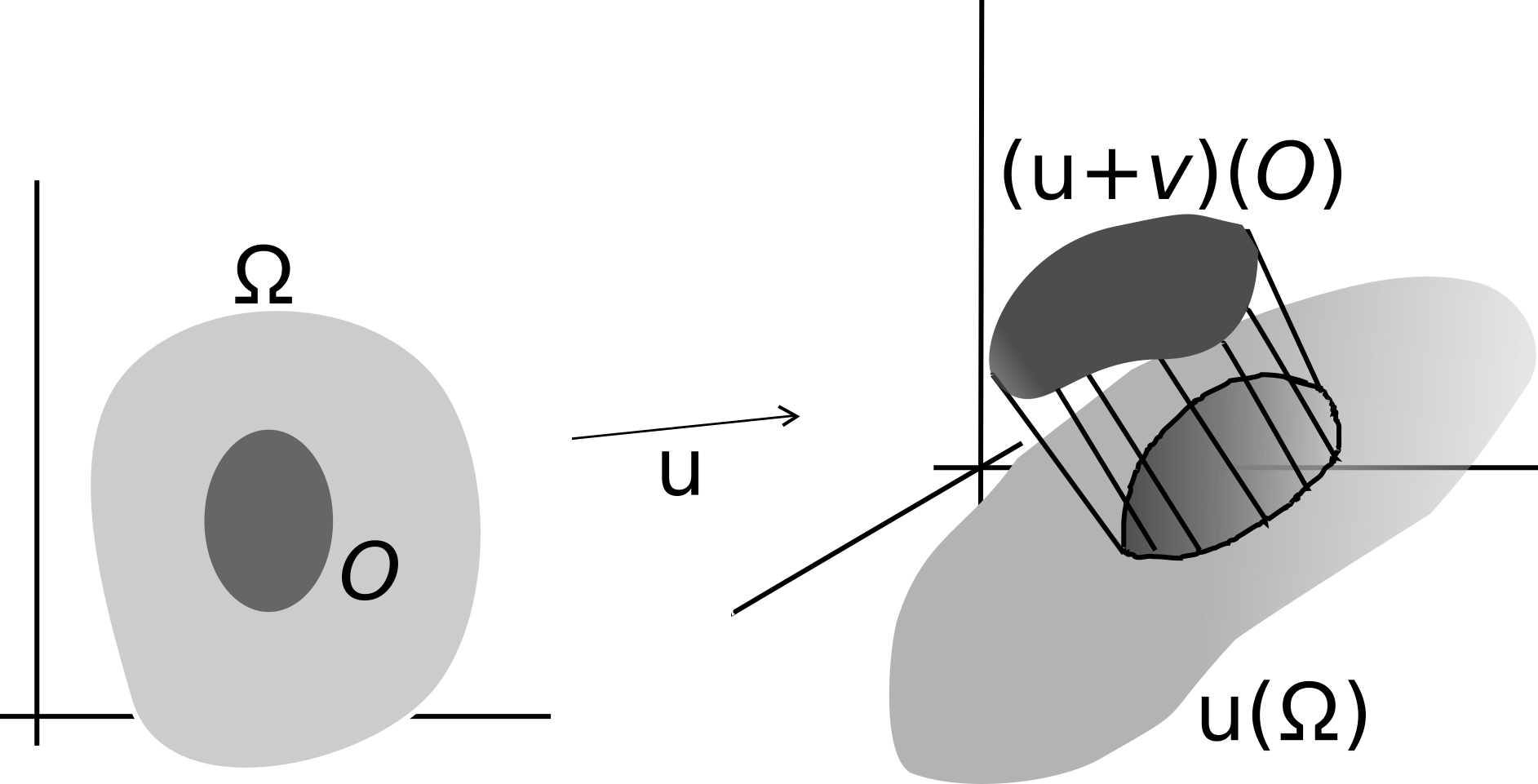}}
\]

\begin{theorem}[Rigidity and flatness of rank-one maps with tangential Laplacian] 
 \label{theorem1}
 
Let $\Om\sub \R^n$ be an open set and $n,N \geq 1$. Let $u \in C^2(\Om,\R^N)$ be a solution to the nonlinear system \eqref{1.1} in $\Om$, satisfying that the rank of its gradient matrix is at most one: 
\[
\rk(\D u) \leq 1 \, \ \text{ in }\Om.
\]
Then, its image $u(\Om)$ is contained in a polygonal line in $\R^N$, consisting of an at most countable union of affine straight line segments (possibly with self-intersections).
\end{theorem}

Let us note that the rank-one assumption for $\D u$ is equivalent to the existence of two vector fields $\xi : \R^n \supseteq \Om \larrow \R^N$ and $a : \R^n \supseteq \Om \larrow \R^n$ such that $\D u= \xi \ot a$ in $\Om$.

Example \ref{example2} below shows that Theorem \ref{theorem1} is optimal and in general rank-one solutions to the system \eqref{1.1} can not have affine image but only piecewise affine.

\begin{example}\label{example2} Consider the $C^2$ rank-one map $u: \R^2 \larrow \R^2$ given by
\[
u(x,y) = \left\{
\begin{array}{rl}
(-x^4,x^4), & x\leq 0,\ y\in \R,
\\
(+x^4,x^4), & x>0,\ y\in \R.
\end{array}
\right.
\]
Then, $u=\nu \circ f$ with $\nu : \R \larrow \R^2$ given by $\nu(t)=(t,|t|)$ and $f : \R^2 \larrow \R$ given by $f(x,y)=\sgn(x)x^4$ (see Figure 2). 
%%%
\[
\underset{\emph{\small Figure 2. The graph of the function $f$ and the image of the curve $\nu$ comprising $u$.}}{\includegraphics[scale=0.18]{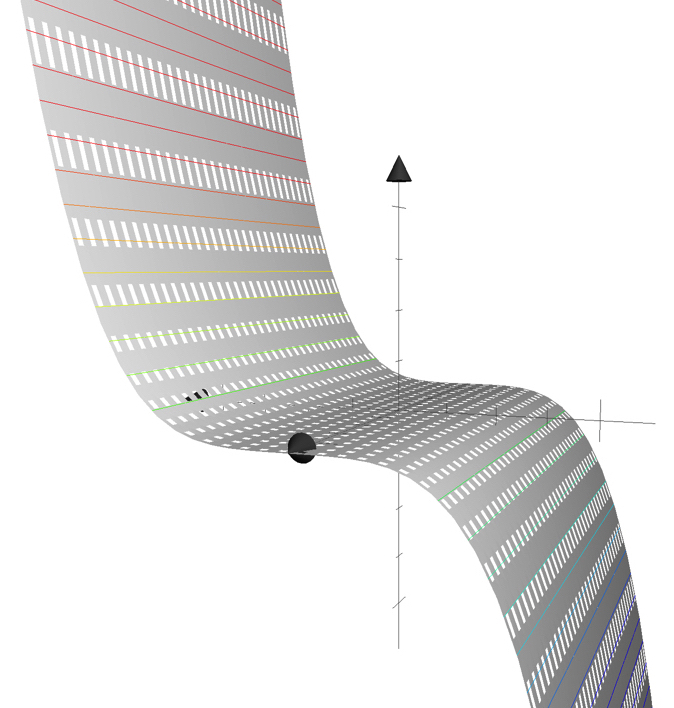} \includegraphics[scale=0.242]{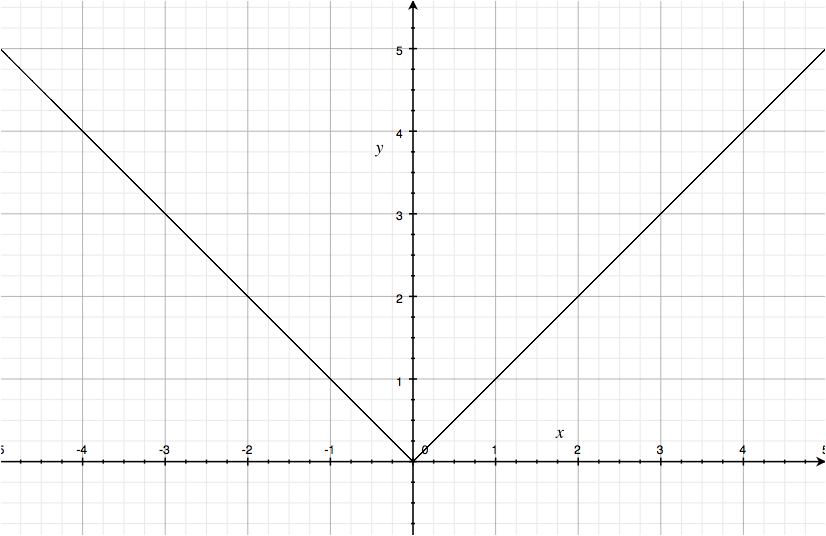}}
\]
It follows that $u$ solves \eqref{1.1} on $\R^2$: indeed, $\De u$ is a non-vanishing vector field on $\{x\neq 0\}$, being tangential to the image thereon since it is parallel to the derivative  $\nu'(t)=(1,\pm 1)$ for $t\neq 0$. On the other hand, on $\{x=0\}$ we have that $\De u =0$. However, the image $u(\R^2)$ of $u$ is piecewise affine but not affine and equals $\nu(\R)$.
\end{example}

As a consequence of Theorem \ref{theorem1}, we obtain the next result regarding the rigidity of $p$-Harmonic maps for $p\in [2,\infty)$ which complements one of the results in the paper \cite{K3} wherein the case $p=\infty$ was considered.

\begin{corollary}[Rigidity of $p$-Harmonic maps, cf.\ \cite{K3}] \label{corollary3} Let $\Om\sub \R^n$ be an open set and $n,N \geq 1$. Let $u \in C^2(\Om,\R^N)$ be a $p$-Harmonic  map in $\Om$ for some $p\in [2,\infty)$, that is $u$ solves \eqref{1.33}. Suppose that the rank of its gradient matrix is at most one: 
\[
\rk(\D u) \leq 1 \, \ \text{ in }\Om.
\]
Then, the same result as in Theorem \ref{theorem1} is true. 

In addition, there exists a partition of $\Om$ to at most countably many Borel sets, where each set of the partition is a non-empty open set with a (perhaps empty) boundary portion, such that, on each of these, $u$ can be represented as 
\[
u = \nu \circ f.
\]
Here, $f$ is a scalar $C^2$ $p$-Harmonic function (for the respective $p\in [2,\infty)$), defined on an open neighbourhood of the Borel set, whilst $\nu : \R \larrow \R^N$ is a Lipschitz curve which is twice differentiable and with unit speed on the image of $f$.
\end{corollary}

Now we move on to discuss our second main result which concerns the rigidity of solutions $u :\R^2 \supseteq \Om \larrow \R^N$ to \eqref{1.1} for $N\geq 2$, having the additively separated form
\beq \label{asp}
u(x,y) \,=\, f(x) + g(y)
\eeq
for some curves $f,g :\R \larrow \R^N$. Solutions of this form are very important in relation to the $\infty$-Laplacian. If $N=1$, all $\infty$-Harmonic functions of this form after a normalisation reduce to the so-called Aronsson solution on $\R^2$
\[
u(x,y) = |x|^{4/3} - |y|^{4/3}
\]
which is the standard explicit example of a non-$C^2$ $\infty$-Harmonic function with conjectured optimal regularity. In the vectorial case, the family of separated solutions is quite large. For $N=2$, a large class of such vectorial solutions was constructed in \cite{K4} and is given by
\[
u(x,y)\, = \int_x^y \! \big(\cos(K(t)),\sin(K(t)) \big) \, \mathrm{d}t
\]
with $K$ a function in $C^1(\R)$ satisfying certain general conditions. The simplest non-trivial example of an $\infty$-Harmonic map with this form (defined on the strip \{$|x-y|<\pi/4 \} \sub \R^2)$ is given by the choice $K(t)=t$. Our second main result asserts that solutions of separated form to \eqref{1.1} have images which are piecewise affine, contained in a union of intersecting planes of $\R^N$. More precisely, we have:

\begin{theorem}[Rigidity and flatness of maps with tangential Laplacian in separated form] 
 \label{theorem2}
 
Let $\Om\sub \R^2$ be an open set and $N \geq 2$. Let $u \in C^{2}( \Om,\R^N)$ be a solution to the nonlinear system \eqref{1.1} in $\Om$, having the separated form $u(x,y)=f(x)+g(y)$, for some $C^{2}$ curves $f, g : \R \larrow \R^N$. 

Then, the image $u(\Om)$ of the solution is contained in an at most countable union of affine planes in $\R^N$. 
\end{theorem}

In addition, the proof of Theorem \ref{theorem2} shows that every connected component of the set $\{\rk(\D u)=2\}$ is contained entirely in an affine plane and every connected component of the set $\{\rk(\D u)\leq 1\}$ is contained entirely in an affine line. 

Note that our result is trivial in the case that $N=n=2$ since the codimension $N-n$ vanishes. Additionally, due to the regularity of the solutions, if a $C^2$ mapping has piecewise affine image, then second derivatives must vanish when first derivatives vanish at the ``breaking points". Further, one might also restrict their attention to domains of rectangular shape, since any map with separated form can be automatically extended to the smallest rectangle containing the domain. 

Also, herein we consider only the illustrative case of $n=2<N$ and do not discuss more general situations, since numerical evidence obtained in \cite{KP} suggests that Theorem \ref{theorem2} does not hold in general for solutions in non-separated form. However, as a consequence of Theorem \ref{theorem2} we have the next particular result:

\begin{corollary}[Rigidity and flatness of $p$-Harmonic maps in separated form] 
 \label{corollary5}
 
Let $\Om\sub \R^n$ be an open set and $n,N \geq 1$. Suppose that $u \in C^2(\Om,\R^N)$ is a $p$-Harmonic  map in $\Om$ for some $p\in [2,\infty]$, that is $u$ solves \eqref{1.33} if $p<\infty$ and \eqref{1.3} if $p=\infty$. 

Then, if $u$ has the separated form $u(x,y)=f(x)+g(y)$ for some $C^{2}$ curves $f,g : \R \larrow \R^N$, the same conclusion as in Theorem \ref{theorem2} holds true.
\end{corollary}

In this paper we try to keep the exposition as simple as possible and therefore we refrain from discussing generalised solutions to \eqref{1.1} and \eqref{1.3} (or \eqref{1.33}). We confine ourselves to merely mentioning that in the scalar case, $\infty$-Harmonic functions are understood in the viscosity sense of Crandall-Ishii-Lions (see e.g.\ \cite{C,K7}), whilst in the vectorial case a new candidate theory for systems has been proposed in \cite{K8} which has already borne significant fruit in \cite{K8,K9, K10, CKP, AyK,KP}.

We now expound on how exactly the nonlinear system \eqref{1.1} arises from \eqref{1.33} and \eqref{1.3}. By expanding the derivatives in \eqref{1.33} and normalising, we arrive at 
\beq \label{1.33A}
 \D u \ot \D u :\D^2u \, +\, \frac{|\D u|^2}{p-2} \De u  = 0. 
\eeq
For any $X\in \R^{N\by n}$, let $[\![X]\!]^\parallel$ denote the orthogonal projection on the range of the linear map $X : \R^n \larrow \R^N$:
\beq \label{1.7}
[\![X]\!]^\|\,:=\, \textrm{Proj}_{\mathrm{R}(X)}.
\eeq
Since the identity of $\R^N$ splits as $\text{I} = [\![\D u]\!]^\| + [\![\D u]\!]^\bot$, by expanding $ \Delta u$ with respect to these projections, 
\[
\D u \ot \D u :\D^2u \, +\, \frac{|\D u|^2}{p-2}[\![\D u]\!]^\|\De u  = -\frac{|\D u|^2}{p-2}[\![\D u]\!]^\bot\De u .
\]
The mutual perpendicularity of the vector fields of the left and right hand side leads via a renormalisation argument (see e.g.\ \cite{K1,K2,K3}) to the equivalence of the $p$-Laplacian with the pair of systems
\beq \label{1.8}
\D u \ot \D u :\D^2u \, +\, \frac{|\D u|^2}{p-2}[\![\D u]\!]^\|\De u  = 0\ , \ \ \  |\D u|^2[\![\D u]\!]^\bot\De u  = 0.
\eeq
The $\infty$-Laplacian corresponds to the limiting case of \eqref{1.8} as $p\to \infty$, which takes the form
\beq \label{1.9}
\D u \ot \D u :\D^2u = 0\ , \ \ \  |\D u|^2[\![\D u]\!]^\bot\De u  = 0.
\eeq
Hence, the $\infty$-Laplacian \eqref{1.3} actually consists of the two independent systems in \eqref{1.9} above. The system $|\D u|^2[\![\D u]\!]^\bot\De u  = 0$ is, at least on $\{\D u \neq 0\}$, equivalent to \eqref{1.1}. Note that in the scalar case of $N=1$ as well as in the case of submersion solutions (for $N\leq n$), the second system trivialises.

We conclude the introduction with a geometric interpretation of the nonlinear system \eqref{1.1}, which can be expressed in a more geometric language as follows:\footnote{This fact has been brought to our attention by Roger Moser.} Suppose that $u(\Omega)$ is a $C^2$ manifold and let $\textbf{A}(u)$ denote its second fundamental form. Then
\[
[\![\D u ]\!]^\perp \Delta u \, = - \, \tr \,\textbf{A}(u)(\D u, \D u). 
\]
The tangential part $[\![\D u ]\!]^\parallel \Delta u$ of the Laplacian is commonly called the \emph{tension field} in the theory of \emph{Harmonic maps} and is symbolised by $\tau(u)$ (see e.g.\ \cite{M}). Hence, we have the orthogonal decomposition 
\[
\Delta u \,=\, \tau(u) \,-\, \tr \, \textbf{A}(u)(\D u, \D u).
\]
Therefore, in the case of higher regularity of the image of $u$, we obtain that the nonlinear system 
\beq \label{eq}
\Delta u \,=\, \tau(u)\ \  \text{ in }\Om,
\eeq
is a further geometric reformulation of our PDE system \eqref{1.1}.

\section{Proofs} \label{section2}

In this section we prove the results of this paper. Before delving into that, we present a result of independent interest in which we represent explicitly the vector field $\A$ arising in the parametric system $\De u =\D u\,\A$, in the illustrative case of $n=2$. 

We will be using the symbolisations ``cof",``det" and ``rk" to denote the cofactor matrix, the determinant function and the rank of a matrix, respectively.

\begin{lemma}[Representation of $\A$]  \label{lemma2}

Let $u \in C^2(\Om,\R^N)$ be given, $\Om\sub \R^2$ open, $N\geq 2$. The following are equivalent:

\begin{enumerate}

\item The map $u$ is a solution to the PDE system \eqref{1.1}.

\smallskip

\item There exists a vector field $\A : \R^2 \supseteq \Om\larrow \R^N$ such that
\[
\De u  = \D u\,\A \ \ \text {in }\Om.
\]
\end{enumerate}

In $\mathrm{(2)}$, as $\A$ one might choose
 \[
\bar{\A} \,:=\, \left\{
\begin{array}{ll}
\displaystyle 
\frac{\mathrm{cof}\big(\D u^\top \D u\big)^{\!\top}}{\mathrm{det}\big(\D u^\top \D u \big)}(\D u)^\top \!\De u, & \text{ on }\{\rk(\D u)=2\},
\ms
\\
\displaystyle (\De u)^\top \! \frac{\D u \, \D u^\top}{|\D u \,\D u^\top|^2}\, \D u, & \text{ on }\{\rk(\D u)=1\},
\\
0, & \text{ on }\{\rk(\D u)=0\}.
\end{array}
\right.
\]

$\A$ is uniquely determined on $\{\rk(\D u)=2\}$ but not on $\{\rk(\D u)<2\}$ and any other $\A$ has the form $\bar{\A}+V$, where $V(x)$ lies in the nullspace of $\D u(x)$, $x\in \Om$. 
\end{lemma}

\BPL \ref{lemma2}. The equivalence between (1)-(2) is immediate, therefore it suffices to show that $\bar \A$ satisfies $\De u =\D u\,\bar \A$ and is unique on $\{\rk(\D u)=2\}$. Let $\A$ be as in (2). On $\{\rk(\D u)=2\}$, the $2\by 2$ matrix-valued map $\D u^\top \D u$ is invertible and
\[
\big(\D u^\top \D u\big)^{-1}\,=\, \frac{\mathrm{cof}\big(\D u^\top \D u\big)^{\!\top}}{\mathrm{det}\big(\D u^\top \D u \big)}.
\]
Since $\D u^\top \!\De u = \D u^\top \D u \,A$, we obtain that $A=\bar A$.

The claim being obvious for $\{\rk(\D u)=0\}=\{\D u=0\}$, it suffices to consider only the set $\{\rk(\D u)=1\}$ in order to conclude. Thereon we have that $\D u$ can be written as 
\[
\ \ \ \D u  = \xi \ot a, \ \ \text{ in }\{\rk(\D u)=1\},
\]
for some non-vanishing vector fields $\xi$ and $a$. By replacing $\xi$ with $\xi|a|$ and $a$ with $a/|a|$, we may assume $|a|\equiv 1$ throughout $\{\rk(\D u)=1\}$. If $\De u =\D u \,\A$, we have $\De u = (\xi \ot a )\A$ and since any component of $\A$ which is orthogonal to $a$ is annihilated, we may replace $\A$ by $\la a$ for some function $\la$. Therefore,
\[
\De u  = (\xi \ot a) \,\A \,=\, (\xi \ot a) (\la a) \,=\, \xi\la|a|^2 = \la \xi 
\]
and hence $\xi \cdot \De u =\la |\xi|^2$ and also $\xi^\top \D u = a |\xi|^2$. On the other hand, since 
\[
\D u \, \D u^\top =\, (\xi \ot a) (a \ot \xi) = \xi \ot \xi, \ \ \ \ \big| \D u \, \D u^\top \big| = |\xi|^2
\]
we infer that 
\[
\A =\, \la a = \left(\frac{\De u \cdot \xi}{|\xi^2|}\right) \left( \frac{\xi^\top \D u }{|\xi^2|}\right)=\, \frac{\De u^\top (\xi \ot \xi)\, \D u}{|\xi \ot \xi|^2}\,=\, (\De u)^\top \! \frac{\D u \, \D u^\top}{|\D u \,\D u^\top|^2}\, \D u,
\]
as claimed. \qed

\ms

We now continue with the proof of the main results. 

\ms

The main analytical tool needed in the proof of Theorem \ref{theorem1} is the next rigidity theorem for maps whose gradient has rank at most one. It was established in \cite{K3} and we recall it below for the convenience of the reader and only in the case needed in this paper.

 \bt[Rigidity of Rank-One maps, cf.\ \cite{K3}] \label{theorem3} Suppose $\Om \sub \R^n$ is an open set and $u$ is in $C^2(\Om,\R^N)$. Then, the following are equivalent:

\ms 

\noi (i) The map $u$ satisfies that $\rk(Du)\leq 1$ on $\Om$. Equivalently, there exist vector fields $\xi : \Om \larrow \R^N$ and $a : \Om \larrow \R^n$ with $a\in C^1(\Om,\R^n)$ and $\xi \in C^1(\Om\set \{a=0\},\R^N)$  such that 
\[
\ \ \D u = \xi \ot a, \ \ \text{ on }\Om.
\] 

\noi (ii) There exists Borel subset $\{B_i\}_{i \in \N}$ of $\Om$ such that 
\[
\Om  = \bigcup_{i=1}^\infty \, B_i
\]
and each $B_i$ equals a non-empty connected open set with a (possibly empty) boundary portion, functions $\{f_i\}_{i\in \N} \in C^2(\Om)$ and curves $\{\nu_i \}_{i \in \N} \sub W^{1,\infty}_{\mathrm{loc}}(\R,\R^N)$ such that, on each $B_i$ the map $u$ has the form
 \beq \label{2.1a}
\ \  u = \nu_i \circ f_i, \ \ \text{ on }B_i.
  \eeq
Moreover, $|\nu_i'|\equiv 1$ on the interval $f_i(B_i)$, $\nu_i'\equiv 0$ on $\R \set f_i(B_i)$ and $\nu_i''$ exists everywhere on $f_i(B_i)$, interpreted as $1$-sided derivative on $\p f_i(B_i)$ (if $f_i(B_i)$ is not open). Also, 
 \beq \label{2.2a}
\ \ \left\{ \ \
 \begin{split}
  \D u \, &=\, (\nu_i' \circ f_i ) \ot \D f_i\ ,  \hspace{122pt} \text{ on }B_i,
  \\
  \D^2 u \, &=\, (\nu_i'' \circ f_i ) \ot \D f_i \ot \D f_i \,+\, (\nu_i' \circ f_i ) \ot \D^2 f_i\ , \ \ \ \text{ on }B_i.
  \end{split}
    \right.
 \eeq
 In addition, the local functions $(f_i)_1^\infty$ extend to a global function $f\in C^2(\Om)$ with the same properties as above if $\Om$ is contractible (namely, homotopically equivalent to a point).
 \et

We may now prove our first main result.

\BPT \ref{theorem1}. Suppose that $u : \R^n \supseteq \Om \larrow \R^N$ is a solution to the nonlinear system \eqref{1.1} in $C^2(\Om,\R^N)$ which in addition satisfies that $\rk(\D u)\leq 1$ in $\Om$. Since $\{\D u =0 \}$ is closed, necessarily its complement in $\Om$ which is $\{\rk(\D u) = 1\}$ is open. 

By invoking Theorem \ref{theorem3}, we have that there exists a partition of the open subset $\{\rk(\D u) = 1\}$ to countably many Borel sets $(B_i)_1^\infty$ with respective functions $(f_i)_1^\infty$ and curves $(\nu_i)_1^\infty$ as in the statement such that \eqref{2.1a}-\eqref{2.2a} hold true and in addition 
\[
\ \ \ \D f_i\, \neq \, 0 \ \ \text{ on }B_i, \ i \in \N. 
\]
Consequently, on each $B_i$ we have
\[
\begin{split}
[\![\D u]\!]^\bot &=\, [\![ (\nu_i' \circ f_i ) \ot \D f_i ]\!]^\bot=\, \mathrm{I} \, -\, \frac{(\nu_i' \circ f_i) \ot (\nu_i' \circ f_i)}{|\nu_i' \circ f_i|^2},
\\
\De u \, &=\, (\nu_i'' \circ f_i ) |\D f_i|^2 \,+\, (\nu_i' \circ f_i ) \De f_i.
\end{split}
\]
Hence, \eqref{1.1} becomes 
\[
\bigg[ \mathrm{I} \, -\, \frac{(\nu_i' \circ f_i) \ot (\nu_i' \circ f_i)}{|\nu_i' \circ f_i|^2}\bigg] \Big( (\nu_i'' \circ f_i ) |\D f_i|^2 \,+\, (\nu_i' \circ f_i ) \De f_i\Big) = 0,
\]
on $B_i$. Since $|\nu_i|^2\equiv 1$ on $f_i(B_i)$, we have that $\nu_i''$ is orthogonal to $\nu_i'$ thereon and therefore the above equation reduces to
\[
\ (\nu_i'' \circ f_i ) |\D f_i|^2  =\, 0 \ \ \text{ on }B_i,\ i \in \N.
\]
Therefore, $\nu_i$ is affine on the interval $f_i(B_i)\sub \R$ and as a result $u (B_i) = \nu_i(f_i(B_i))$ is contained in an affine line of $\R^N$, for each $i\in \N$. On the other hand, since
\[
u(\Om)\, = \, u\big(\{\D u =0\} \big) \bigcup_{i\in \N}\, u(B_i)
\]
and $u$ is constant on each connected component of the interior of $\{\D u =0\}$, the conclusion ensues by the regularity of $u$ because $u\big(\{\D u =0\} \big)$ is also contained in the previous union of affine lines. The result ensues. \qed

\ms

Now we establish Corollary \ref{corollary3} by following similar lines to those of the respective result in \cite{K3}.

\BPCOR \ref{corollary3}. Suppose $u$ is as in the statement of the corollary. By Theorem \ref{theorem3}, there exists, a partition of $\Om$ to Borel sets $\{B_i\}_{i \in \N}$, functions $f_i\in C^2(\Om)$ and Lipschitz curves $\{\nu_i\}_{i \in \N} : \R \larrow \R^N$ with $|\nu_i'|\equiv 1$ on $f_i(B_i)$, $|\nu_i'|\equiv 0$ on $\R \set f_i(B_i)$ and twice differentiable on $f_i(B_i)$, such that $u|_{B_i} = \nu_i \circ f_i$ and \eqref{2.2a} holds as well. Since on each $B_i$ we have 
 \[
 |\D u|  = \big|(\nu_i' \circ f_i) \ot \D f_i \big| = |\D f_i|,
\]
by \eqref{1.33A} and the above, we obtain 
\[
\begin{split}
\big( (\nu_i' \circ f_i) \ot \D f_i & \big) \ot  \big( (\nu_i' \circ f_i) \ot \D f\big) : \bigg[ (\nu_i'' \circ f_i) \ot \D f_i \ot \D f_i \, +\, (\nu_i' \circ f_i) \ot \D^2f_i
\bigg] 
 \\
 & + \, \frac{|\D f_i|^2}{p-2}\bigg\{ (\nu_i' \circ f_i)\, \De f_i\,+\, (\nu_i'' \circ f) |\D f_i|^2\bigg\}\ =\ 0,
\end{split}
\]
on $B_i$. Since $\nu_i''$ is orthogonal to $\nu'_i$ and also $\nu'_i$ has unit length, the above reduces to 
\[
\begin{split}
(\nu_i' \circ f_i) \bigg[  \D f_i  \ot Df_i  : \D^2f_i \, +\, \frac{|\D f_i|^2}{p-2}\, \De f_i \bigg] 
 \, + \, \frac{1}{p-2}(\nu_i'' \circ f_i) |\D f_i|^4 \ =\ 0,
\end{split}
\]
on $B_i$. Again by orthogonality, the above is equivalent to the pair of independent systems
\[
\begin{split}
(\nu_i' \circ f_i) \bigg[  \D f_i  \ot Df_i  : \D^2f_i \, +\, \frac{|\D f_i|^2}{p-2}\, \De f_i \bigg] 
 \, =\,0\ , \ \ \  (\nu_i'' \circ f_i) |\D f_i|^4 \ =\ 0,
\end{split}
\]
on $B_i$. Since $|\nu_i'|\equiv 1$ of $f_i(B_i)$, it follows that $\De_p f_i = 0$ on $B_i$ and since $(B_i)_1^\infty$ is a partition of $\Om$ of the form described in the statement, the result ensues by invoking Theorem \ref{theorem1}.          \qed

\ms

We may now prove our second main result.

\BPT \ref{theorem2}. We begin by temporarily assuming that $\Om$ is a rectangle of the form
\[
Q = (a,b)\by(c,d)\, \sub\, \R^2
\]
and we fix $(x_0,y_0)\in Q= \Om$. Let us also assume that the rank of the gradient is full throughout:
\[
\rk(\D u)\, \equiv \, 2 \ \ \text{ in }Q.
\]
Later we will remove both these extra assumptions. Let $f:(a,b) \longrightarrow \R^N$ and $g:(c,d) \longrightarrow \R^N$ be such that $u(x,y)=f(x)+g(y)$. Then, the gradient matrix then has the form 
\[
\D u(x,y)  = \big( f'(x) , g'(y)\big) \, \in \, \R^{N \times 2}.
\]
By assumption, we have that  $[\![\D u]\!]^\bot \Delta u= 0$ in $\Om$. By Lemma \ref{lemma2}, there exists  a vector field $\A: \R^2 \supseteq \Om \longrightarrow \R^2$ such that $\Delta u =\D u \, \A$. If $\A$ has components $(a,b)$, this means that the functions $f$ and $g$ satisfy
\beq \label{2.2}
f''(x) \,+\, g''(y) \,=\, a(x,y) f'(x) \,+\, b(x,y) g'(y).
\eeq
Although we will not utilise this in the sequel, it is instructive and quite possible to express the coefficients $a,b$ in terms of $f,g,f',g'$ along the lines of Lemma \ref{lemma2} but more concretely, as follows. By applying the $\R^{N\by N}$ matrix
\[
|g'(y)|^2\,\mathrm{I}-g'(y)\ot g'(y) 
\]
to \eqref{2.2}, the summand $b(x,y) g'(y)$ on the right hand side is annihilated and we obtain
\[
\Big[|g'(y)|^2\, \mathrm{I}-g'(y)\ot g'(y)\Big]\big(f''(x) \,+\, g''(y)\big) = a(x,y) \Big[|g'(y)|^2\, \mathrm{I}-g'(y)\ot g'(y)\Big] f'(x).
\]
Hence, on $\{\rk(\D u)=2\}$ we have
\[
a(x,y)  = \frac{|g'(y)|^2\, \mathrm{I}-g'(y)\ot g'(y)}{\big[|g'(y)|^2\, \mathrm{I}-g'(y)\ot g'(y)\big] f'(x)}\big(f''(x) \,+\, g''(y)\big).
\]
Arguing symmetrically, we may obtain
\[
b(x,y)  = \frac{|f'(x)|^2\, \mathrm{I}-f'(x)\ot f'(x)}{\big[|f'(x)|^2\, \mathrm{I}-f'(x)\ot f'(x)\big] g'(y)}\big(f''(x) \,+\, g''(y)\big).
\]
Similar expression can be obtained on $\{\rk(\D u)\leq 1\}$ as well.

Let us consider \eqref{2.2} as a first order ODE with unknown function $f'$ in the variable $x$. By integrating the equation in $x$, we view its first integral as an first order ODE with solution the function $g'$ in the variable $y$, which can be integrated again. Therefore, by sparing the reader  the tedious computations, we arrive at the following integral identity
\beq \label{2.3}
\left\{\ \ 
\begin{split} 
 g'(y) \, e^{- \int_{y_0}^y \frac{B(x,s)}{A(x,s)}\,\mathrm{d}s} \, & + \, f'(x)  \int_{y_0}^y \frac{ e^{- \int_{x_0}^x a(s,t)\,\mathrm{d}s}}{A(x,t)} \, e^{- \int_{y_0}^t \frac{B(x,\tau)}{A(x,\tau)}\,\mathrm{d}\tau}\,\mathrm{d}t 
 \\
&=\, g'(y_0) \,+\, f'(x_0)  \int_{y_0}^y \frac{1}{A(x,t)} e^{- \int_{y_0}^t \frac{B(x,\tau)}{A(x,\tau)}\,\mathrm{d}\tau}\,\mathrm{d}t,
\end{split}
\right.
\eeq
where 
\[
\left\{ \ \ 
\begin{split} 
A(x,y)\, &:=\int_{x_0}^x e^{- \int_{x_0}^s a(\tau,y)\,\mathrm{d}\tau}\,\mathrm{d}s\ , 
\\ 
B(x,y)\, & := \int_{x_0}^x b(s,y) \, e^{- \int_{x_0}^s a(\tau,y)\,\mathrm{d}\tau}\,\mathrm{d}s.
\end{split}
\right.
\]
Note that 
\[
\text{$A(x,y)\lessgtr 0$ \ \ if and only if \ \ $x-x_0\lessgtr 0$, \ for all $y\in \R$.}
\]
Let us now define for convenience the functions
\[
\left\{\ \ \
\begin{split}
C (x,y )\, &:=\,e^{- \int_{y_0}^y \frac{B(x,s)}{A(x,s)}\,\mathrm{d}s},
\\
D  (x,y )\, & := \int_{y_0}^y \frac{ e^{- \int_{x_0}^x a(s,t)\,\mathrm{d}s}}{A(x,t)} e^{- \int_{y_0}^t \frac{B(x,\tau)}{A(x,\tau)}\,\mathrm{d}\tau} \,\mathrm{d}t, 
\\
E  (x,y )\, &:= \int_{y_0}^y \frac{1}{A(x,t)} e^{- \int_{y_0}^t \frac{B(x,\tau)}{A(x,\tau)}\,\mathrm{d}\tau} \,\mathrm{d}t.
\end{split}
\right.
\]
For the sake of brevity, we will be suppressing the dependence on the variables $(x,y )$. In view of these definitions, \eqref{2.3} can be re-written in the abbreviated form 
\beq \label{2.4}
g'(y)  C \, +\, f'(x)  D  \,= \, g'(y_0) \, +\, f'(x_0) E .
\eeq
By the symmetry of the equations, the substitutions 
\[
\text{ $g  \longleftrightarrow f,\ x  \longleftrightarrow y,\ x_0  \longleftrightarrow y_0,\ a  \longleftrightarrow b$ }
\]
lead to the additional equation
\beq  \label{2.5}
\left\{ \ \ 
\begin{split} 
f'(x) \, e^{- \int_{x_0}^x \frac{G (s,y)}{F (s,y)}\,\mathrm{d}s} \, &+ \, g'(y) \int_{x_0}^x \frac{ e^{- \int_{y_0}^y b(t,s)\,\mathrm{d}s}}{F (t,y)}  e^{- \int_{x_0}^t \frac{G (\tau,y)}{F (\tau,y)}\,\mathrm{d}\tau}\,\mathrm{d}t 
\\ 
&=\, f'(x_0)\, +\, g'(y_0)  \int_{x_0}^x \frac{1}{F (t,y)} e^{- \int_{x_0}^t \frac{G (\tau,y)}{F (\tau,y)}\,\mathrm{d}\tau}\,\mathrm{d}t,
\end{split}
\right.
\eeq
where
\[
\left\{ \ \ 
\begin{split}
F (x,y)\, &:=\, \int_{y_0}^y e^{- \int_{y_0}^s b(x,\tau)\,\mathrm{d}\tau}\,\mathrm{d}s\ ,
\\
G (x,y)\,&:=\, \int_{y_0}^y a(x,s)  e^{- \int_{y_0}^s b(x,\tau)\,\mathrm{d}\tau}\,\mathrm{d}s. 
\end{split}
\right. 
\]
Note that
\[
\text{ $F (x,y)\lessgtr 0$ \ \ if and only if \  \ $y-y_0 \lessgtr0$ \ for all $x\in \R$.}
\]
For brevity, we set 
\[
\left\{ \ \ 
\begin{split}
H   (x,y )\, & :=\,e^{- \int_{x_0}^x \frac{G (s,y)}{F (s,y)}\,\mathrm{d}s},
\\
I  (x,y )\, &:=\, \int_{x_0}^x \frac{ e^{- \int_{y_0}^y b(t,s)\,\mathrm{d}s}}{F (t,y)} e^{- \int_{x_0}^t \frac{G  (\tau,y)}{F (\tau,y)}\,\mathrm{d}\tau} \,\mathrm{d}t,\nonumber
\\
J  (x,y )\,&:=\, \int_{x_0}^x \frac{1}{F (t,y)} e^{- \int_{x_0}^t \frac{G  (\tau,y)}{F (\tau,y)}\,\mathrm{d}\tau} \,\mathrm{d}t.
\end{split}
\right.
\]
Then  \eqref{2.5} can be re-written in the simpler form 
\beq  
f'(x) \,=\, {H  }^{-1}\Big( g'(y_0)  J \, +\, f'(x_0) \,-\, g'(y)   I  \Big).\nonumber
\eeq
Substituting the above into \eqref{2.4}, after some elementary calculations we obtain the equation
\beq \label{2.7}
\big(C -  {I } {H  }^{-1}  \big) \,  g'(y)  = \big( E  -  {D } {H  }^{-1}  \big) f'(x_0) \, + \, \big( 1+ {J } {H  }^{-1}   \big) g'(y_0), 
\eeq
for all $(x,y)\in Q$. Note that 
\[
I  < 0 \ \  \text{ if }\ \ (x-x_0)(y-y_0) < 0,\ \  (x,y)\in Q. 
\]
Moreover, we have that $C > 0$ and also $ H   > 0$, for all $(x,y)\in Q$. Hence, if $y < y_0$, we may choose $x > x_0$ and if $y  > y_0$, we may choose $x< x_0$. In either case, we can arrange 
\[
C - {I }{H  }^{-1} \,>\, 0
\]
for all $ y \ne y_0$ such that $(x-x_0)(y-y_0) < 0$. Therefore, from \eqref{2.7} we deduce that 
\[
  g'(y)  = \left(\frac{E-{D } {H  }^{-1}}{C -  {I } {H  }^{-1}}\right) f'(x_0) \, + \, \left(\frac{ 1+ {J } {H  }^{-1} }{ C -  {I } {H  }^{-1}}\right) g'(y_0), 
\]
which yields
\[
g'(y) \,\in\, \text{span} [ \big\{ f'(x_0),g'(y_0)\big\} ]. 
\]
By an integration, the above inclusion implies that the curve $g$ is valued in an affine plane of the form
\[
\phantom{\Big|} g(y) \, \in \,  g(y_0)+\text{span} [\big\{ f'(x_0),g'(y_0) \big\} ].
\]
By arguing similarly for $f$, we also infer that 
\[
 f'(x)\, \in \, \text{span} [ \big\{ f'(x_0),g'(y_0) \big\} ]
\] 
and hence
\[
\phantom{\Big|}  f(x) \, \in \,  f(x_0)+ \text{span} [\big\{ f'(x_0),g'(y_0)\big\} ].
\]
Conclusively, by putting the above together we have obtained
\[
  f(x)+g(y) \, \in \,  \big(f(x_0)+g(y_0)\big)+\text{span} [ \big\{ f'(x_0),g'(y_0) \big\} ],
\]
that is,
\[
\phantom{\Big|} u(x,y) \, \in \, u(x_0,y_0) +  \mathrm{R}\big(\D u(x_0,y_0)\big),
\]
for all $(x,y)$ running in the rectangle $Q$. Hence, $u(Q)$ is contained in an affine plane of $\R^N$, as claimed.

\smallskip

Finally, we remove the additional assumptions that $\Om$ is a rectangle and that the rank of $\D u$ is full. Let us consider separately each connected component of the open subset of $\Om$ given by $\{\rk(\D u)=2\}$. We cover each component with countably many (overlapping) rectangles $\big\{Q_i : i\in \N\big\}$ where 
\[
Q_i  = (a_i, b_i)\by(c_i, d_i), \ \ \ i\in \N ,
\]
and respective points inside the rectangles 
\[
\big\{(x_{0i},y_{0i}) \in Q_i \ :  \ i \in \N \big\}. 
\]
On each rectangle, by the previous reasoning the solution will be contained in an affine plane $\Pi_i \sub \R^N$. By choosing $(x_{0i},y_{0i})$ on the overlaps of each $Q_i$ with all its neighbouring rectangles, connectedness of the component allows us to conclude that all the planes coincide. 

It remains to consider the complement $\Om\set \{\rk(\D u)=2\}$, which we decompose to the sets
\[
\inter  \big\{\rk(\D u)\leq 1 \big\}\, \bigcup \, \p \big\{\rk(\D u)\leq 1 \big\},
\]
where ``int" denotes topological interior. If the interior is non-empty, on each connected component of it we apply Theorem \ref{theorem1} to infer that $u\big( \inter \{\rk(\D u)\leq 1 \}\big)$ is contained in a polygonal curve of $\R^N$, given by an at most countable union of affine straight line segments. Finally, by continuity of $u$, we have that $u\big( \p \{\rk(\D u)\leq 1 \}\big)$ is also contained in the union of planes. The result ensues.      \qed
 
 \ms

We conclude by establishing the remaining corollary.

\BPCOR \ref{corollary5}. Suppose that $u$ is a $C^2$ $p$-Harmonic mapping as in the hypotheses of the corollary for some $p\in[2,\infty]$, namely that it has additively separated form (see \eqref{asp}) and solves either \eqref{1.33} if $p<\infty$ or \eqref{1.3} if $p=\infty$. By \eqref{1.8}-\eqref{1.9}, we deduce that $u$ solves the PDE system
\[
\ |\D u|^2[\![\D u]\!]^\bot\De u  = 0 \ \ \text{ in } \Om. 
\]
On the open set $\{\D u \neq 0\}$, we readily have that $u$ solve the system \eqref{1.1}. On the other hand, we decompose its complement to
\[
\inter  \big\{\D u=0\big\}\, \bigcup \, \p \big\{\D u=0\big\},
\]
If $\inter \{\D u=0 \}$ is non-empty, then $u$ is constant on each connected component of it. Finally, again by the regularity of $u$, we have that $u\big( \p \{\D u=0 \}\big)$ is also contained in the previous union of planes. The claim has been established.     \qed
\ms

\subsection*{Acknowledgement.} N.K.\ is indebted to G. Aronsson for his encouragement and support towards him. He would also like to thank Roger Moser and Giles Shaw for their comments and suggestions which improved the content of this paper. We would also like to thank Jan Kristensen, Tristan Pryer and Enea Parini for several scientific discussions relevant to the results herein.

\bibliographystyle{amsplain}

\ms

\end{document}